\documentclass{amsart}
\usepackage{amscd}
\usepackage{graphicx}
\newcommand{\bbR}{{\mathbb R}}
\newcommand{\bbZ}{{\mathbb Z}}
\newtheorem{lemma}{Lemma}[section]
\newtheorem{prop}[lemma]{Proposition}
\newtheorem{thm}{Theorem}

\begin{document}
\title{Minimal spheres of arbitrarily high Morse index}
\author{Joel Hass}
\address{Department of Mathematics\\University of California, Davis\\
CA 95616\\USA}
\email{hass@math.ucdavis.edu}
\thanks{PN and JHR would like to acknowledge the support of the 
Australian Research Council}
\author{Paul Norbury}
\address{Department of Mathematics and Statistics\\
University of Melbourne\\Australia 3010}
\email{pnorbury@ms.unimelb.edu.au}
\author{J. Hyam Rubinstein}
\address{Department of Mathematics and Statistics\\
University of Melbourne\\Australia 3010}
\email{rubin@ms.unimelb.edu.au}

\begin{abstract}
We construct a smooth Riemannian metric on any 3-manifold with the
property that there are genus zero embedded minimal surfaces of
arbitrarily high Morse index.
\end{abstract}

\maketitle

\section{Introduction.}
In connection with the spherical space form problem, Pitts and
Rubinstein \cite{PRuApp} raised the question of whether embedded
minimal surfaces of fixed genus in a three-manifold have bounded Morse
index.  An area bound on minimal surfaces of fixed genus would imply
an index bound, however in \cite{CMiExa} it was shown that on any
three-manifold there exists an open set of metrics for which there are
embedded minimal tori of arbitrarily large area.  These minimal tori
have Morse index zero.  In this paper we construct metrics for which
there are embedded minimal spheres of arbitrarily large Morse index.

We use the techniques from \cite{HLaMin} to study circle invariant
metrics and reduce the problem to one of finding a family of geodesics
of arbitrarily large index on a surface equipped with a metric that
degenerates at the boundary.

\begin{thm}
    On any three-manifold there exists a metric for which there are
    embedded minimal spheres of arbitrarily large Morse index.  On the
    three-sphere it can be chosen to have non-negative Ricci
    curvature.
\end{thm}

The Ricci curvature of each example of the theorem is non-positive
somewhere.  Each of the metrics on the three-sphere constructed by the
theorem possesses an arbitrarily small perturbation to a positive
Ricci curvature metric.  If a manifold has a metric of positive Ricci
curvature then there is an area bound on the embedded minimal spheres
\cite{CScSpa,CMiMin} and hence a bound on their index.  Thus the
property of the existence of embedded minimal spheres of arbitrarily
large Morse index is destroyed by the perturbations.

Using completely different methods, Colding and Hingston \cite{CHiMet}
have constructed metrics on any three-manifold with embedded minimal
tori of arbitrarily high index.  This is complementary to the work
here since the conjecture of Pitts and Rubinstein \cite{PRuApp}
requires a Morse index bound on a compactification of the space of
embedded tori which includes both spheres and tori.

\section{The quotient construction.}
Consider the product metric on $\bbR\times S^2$ given by the Euclidean 
line times the round two-sphere: 
\[ds_0^2=dr^2+d\phi^2+{\rm sin}^2\phi\hspace{2 pt}d\theta^2\] 
where $\phi\in[0,\pi]$ and $\theta\in[0,2\pi]$ are spherical
coordinates on the two-sphere.  This is invariant under rotations by 
$\theta$, and the quotient metric is 
\[ds_1^2=dr^2+d\phi^2\]
defined on $\bbR\times[0,\pi]$.  Geodesics on the quotient surface 
have little to do with minimal surfaces upstairs since the lengths of 
the orbit circles vary.  Instead, following \cite{HLaMin} we adjust 
the metric by $\sin\phi$, the lengths of the orbit circles:
\begin{equation}  \label{eq:degmet}
    ds^2=\sin^2\phi(dr^2+d\phi^2).
\end{equation}
This is a metric on $\bbR\times[0,\pi]$ degenerate on the boundary, 
and its geodesics pull back to minimal surfaces in $\bbR\times S^2$.
Its curvature is $1/\sin^4\phi$.

\begin{figure}[ht]
\begin{center}
\scalebox{0.5}{\includegraphics[height=8cm]{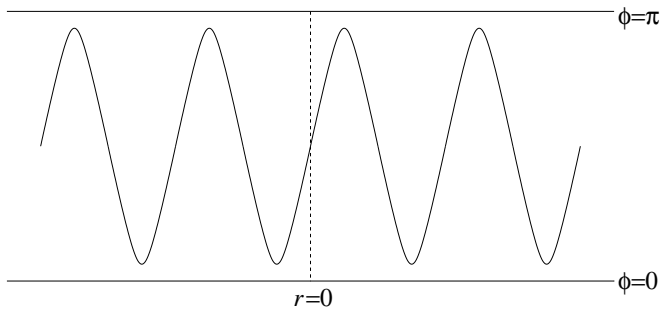}}
\end{center}
\caption{The degenerate boundary repels geodesics.}
\label{fig:repel}
\end{figure}

The equation for geodesics is given by:
\begin{eqnarray} \label{eq:geod1}
    \ddot{r}&=&-2\frac{\cos\phi}{\sin\phi}\dot{r}\dot{\phi}\\
    \ddot{\phi}&=& 
    \frac{\cos\phi}{\sin\phi}\left(\dot{r}^2-\dot{\phi}^2\right)
\end{eqnarray}
and (\ref{eq:geod1}) is equivalent to 
$d(\ln\dot{r})/dt=-2d(\ln\sin\phi)/dt$ so 
\begin{equation} \label{eq:rdot}
    \dot{r}=\frac{c}{\sin^2\phi}
\end{equation}
for a constant $c$.  The constant $c$ is equal to zero precisely when
the geodesic is vertical ($r=$ constant) and thus meets the boundary. 
For all other geodesics $\dot{r}$ is always positive or always
negative.  The speed $(\dot{r}^2+\dot{\phi}^2)\sin^2\phi$ is a
conserved quantity and we may set it to be $1$.  Thus, $|\dot{r}|\leq
1/\sin\phi$ so from (\ref{eq:rdot}),
\begin{equation}  \label{eq:repel}
    0\leq|c|\leq\sin\phi.
\end{equation}
For non-vertical geodesics $|c|$ is strictly positive so $\phi$ is
bounded away from $0$ and $\pi$, or in other words the boundary repels
the geodesic---see Figure~\ref{fig:repel}.  The maximum and minimum
values of $\phi$ of the geodesic satisfy $\sin\phi=|c|$ and the
geodesic is periodic by the translation and reflection invariance of
$ds^2$.

This is much like the construction in \cite{HLaMin} for the round
metric on $S^3$.  The major difference is that the period of the
geodesics can be made arbitrarily small in our case whereas there is a
lower bound for the period in \cite{HLaMin}.
\begin{lemma}  \label{th:period}
    The period, in $r$, of the geodesic with minimum $\sin\phi=c>0$ is
    less than $2\pi\sqrt{2c/(1+c)}$.
\end{lemma}
\begin{proof}
    From (\ref{eq:rdot}) and $(\dot{r}^2+\dot{\phi}^2)\sin^2\phi=1$ we 
    get
    \[\left(\frac{dr}{d\phi}\right)^2=\frac{\dot{r}^2}{\dot{\phi}^2}
    =\frac{c^2}{\sin^2\phi-c^2}.\]
    By symmetry one quarter of a period of the geodesic is given by 
    the following integral:
    \begin{eqnarray*}
	r&=&\int^{\pi/2}_{\sin\phi=c}\frac{c}{\sqrt{\sin^2\phi-c^2}}d\phi\\
	&=&c\int^1_c\frac{du}{\sqrt{(u^2-c^2)(1-u^2)}}\\
	&<&\sqrt{\frac{c}{2(1+c)}}\int^1_c\frac{du}{\sqrt{(u-c)(1-u)}}\\
	&=&\pi\sqrt{\frac{c}{2(1+c)}}.
    \end{eqnarray*}
\end{proof}
The {\em energy functional} on the space of paths between points $p$
and $q$ in a manifold $M$ is defined to be the integral of the square
of the norm of the derivative along the path.  If we simply integrate
the norm along the path then this gives the length functional.  The
critical points of the energy functional coincide with the critical
points of the length functional, and thus are the geodesics.  The 
Hessian of the energy functional is a bilinear form with 
finite index and nullity---the dimension of the largest negative 
eigenspace, respectively the null space.  It is non-degenerate if the 
nullity is zero.

For generic $p$ and $q$ the energy functional is Morse.  That is, its
critical points are non-degenerate \cite{MilMor}.  Along any geodesic
$\gamma(t)$, the energy functional can be defined on an initial
interval $t\in[0,\tau]$ say.  The geodesic is non-degenerate for
generic $\tau$, and its index remains the same for nearby generic
values of $\tau$.  It is important to note that the energy functional
is defined on a different set of paths for each $\tau$ and in
particular that the tangent spaces change---the infinitesimal
variations are given by vector fields that vanish at $0$ and $\tau$. 
On each side of a non-generic value of $\tau$ the index of the
geodesics can be related, leading to the Morse index theorem.

The points $p$ and $q$ are conjugate along a geodesic $\gamma$ joining
them if the nullity of $\gamma$ is non-zero, or in other words the
geodesic is a degenerate critical point of the energy functional. 
Thus along $\gamma$ there exists a null vector field or equivalently a
{\em Jacobi vector field} that vanishes at $p$ and $q$.  A Jacobi
vector field is a vector field $W$ defined along $\gamma$ that
satisfies the Jacobi equation: \[ D^2W/dt^2+R(W,V)V=0\] where $D$ is
the covariant derivative along $\gamma$, $V=d\gamma/dt$ and
$R(u,v)=D_vD_u-D_uD_v$ is the curvature of the metric.  More
generally, any two points on $\gamma$ are conjugate along $\gamma$ if
there exists a Jacobi vector field along $\gamma$ that vanishes at the
two points.

The Morse index theorem states that the index of a geodesic is equal 
to the number of points conjugate to one of its endpoints.  The 
conjugate points are counted with multiplicity which is $1$ on a 
surface.

On the surface with metric (\ref{eq:degmet}) there are geodesics with
arbitrarily high index.  This uses Lemma~\ref{th:period} which tells
us that for any $N\in\bbZ$, when $c$ is small enough the geodesic with
minimum $\sin\phi=c$ meets $\phi=\pi/2$ at least $N$ times, and the
next lemma in which we show that this gives a lower bound of $N$ for
the Morse index of the geodesic.
\begin{lemma}
    The index of a geodesic is bounded below by half the number of
    times it meets $\phi=\pi/2$.
\end{lemma}
\begin{proof}
    Along any geodesic $\gamma$ with minimum $\sin\phi=c$, choose a
    piece $\gamma^0$ that travels a full period.  In the next
    paragraph we show there exists a pair of conjugate points along
    $\gamma^0$, and thus by the Morse index theorem, the index of
    $\gamma^0$ is at least 1.  That is, there exists a vector field
    along $\gamma^0$ that vanishes at its endpoints and on which the
    Hessian of the energy functional is negative definite.  Extend
    this to $\gamma$ by setting it to be zero outside $\gamma^0$.  Do
    this along consecutive pieces of the geodesic $\gamma$ to get a
    set of vector fields on which the Hessian of the energy functional
    is diagonal and negative definite.  Hence the index of a geodesic
    is bounded below by half the number of times it meets
    $\phi=\pi/2$.
    
    Since the Gaussian curvature of the metric is bounded below by 1
    we can use the Rauch comparison theorem to deduce that the
    distance between conjugate points is at most $\pi$.
    The length of $\gamma^0$ is given by:
    \begin{eqnarray*}
	l&=&2\int^{\sin\phi=c}_{\phi=\pi/2}\sin\phi
	\sqrt{1+(dr/d\phi)^2}\hspace{2 pt}d\phi\\
	&=&2\int^{\sin\phi=c}_{\phi=\pi/2}
	\frac{\sin^2\phi}{\sin^2\phi-c^2}d\phi\\
	&=&2\int^{\pi/2}_0\sqrt{\cos^2\alpha+c^2\sin^2\alpha}\hspace{2
        pt}d\alpha\\
	&\geq&2\int^{\pi/2}_0\cos\alpha\hspace{2 pt}d\alpha=2
    \end{eqnarray*}
    where we have used the substitution $\sin\phi=\sin\alpha\hspace{2
    pt}\sqrt{1-c^2}$.  Thus the length of a geodesic along its full
    period is greater than 4 so contains conjugate points by the Rauch
    comparison theorem.
\end{proof}
We can improve the previous lemma and bound the index of a geodesic
below by the number of times it meets $\phi=\pi/2$.  A minimax
argument shows that the index of a piece of the geodesic that travels
from $\phi=\pi/2$ to its minimum and back---half the period---is 1. 
This uses a sequence of sweepouts between two index zero
geodesics---the horizontal geodesic along $\phi=\pi/2$ and the
degenerate piecewise geodesic consisting of two vertical geodesics
from $\phi=\pi/2$ to $\phi=0$ together with a path (of zero length)
along the boundary.  The minimax is an index 1 geodesic.

\section{Perturbed metric.}
Perturb the product metric on $\bbR\times S^2$ to get
\[ds_0^2=dr^2+\lambda(r)(d\phi^2+{\rm sin}^2\phi d\theta^2)\] where
$\lambda(r)=1$ for $r\leq 0$ and $\lambda'(r)<0$ for $r>0$.  The
metric is defined only for $r$ satisfying $\lambda(r)\geq 0$.  We also
require that $\lambda'(r)=0$ when $\lambda(r)=0$ so that the metric 
is smooth as $\lambda(r)\rightarrow 0$.

The metric on $\bbR\times[0,\pi]$ degenerate on the boundary is
\begin{equation}
    ds^2=\lambda(r)\sin^2\phi(dr^2+\lambda(r)d\phi^2).
\end{equation}
The equation for geodesics is given by:
\begin{eqnarray} \label{eq:geod2}
    \ddot{r}&=&-2\frac{\cos\phi}{\sin\phi}\dot{r}\dot{\phi}
    -\frac{\lambda'(r)}{2\lambda(r)}\dot{r}^2
    +\lambda'(r)\dot{\phi}^2\\
    \ddot{\phi}&=&\frac{\cos\phi}{\sin\phi}\left(\frac{1}{\lambda(r)}\dot{r}^2
    -\dot{\phi}^2\right)-2\frac{\lambda'(r)}{\lambda(r)}\dot{r}\dot{\phi}
\end{eqnarray}
or, alternatively
\begin{equation}  \label{eq:implicit}
    \frac{d^2r}{d\phi^2}=-\frac{1}{\lambda(r)}\frac{\cos\phi}{\sin\phi}
    \left(\frac{dr}{d\phi}\right)^3
    +\frac{3}{2}\frac{\lambda'(r)}{\lambda(r)}\left(\frac{dr}{d\phi}\right)^2
    -\frac{\cos\phi}{\sin\phi}\frac{dr}{d\phi}+\lambda'(r).
\end{equation}
Put
\begin{equation}  \label{eq:c1}
    \lambda(r)=\cos^2r
\end{equation}
for $0\leq r\leq\pi/2$, and $\lambda\equiv 1$ for $r<0$.  This gives a
$C^1$ metric.  A family of solutions to (\ref{eq:implicit}) is
\begin{equation}  \label{eq:explicit}
    \cos\phi=\frac{\tan r}{\tan\kappa}
\end{equation}
for $r\geq 0$ and a constant $\kappa\in[0,\pi]$.  We uniquely continue
to a periodic solution on $r<0$.  Each of these geodesics meets the
boundary at right angles, at $(\kappa,0)$ or $(\pi-\kappa,\pi)$, and
passes through the focal point $(r,\phi)=(0,\pi/2)$ at an angle of
$\theta(\kappa)$ with the line $r=0$.  As $\kappa\rightarrow 0$,
$\theta(\kappa)\rightarrow 0$.

For a given $\epsilon>0$, countably many of these solutions meet
$(r,\phi)=(-\epsilon,\pi/2)$ since as $\theta(\kappa)$ gets smaller,
so does the period and geodesics with a period of $2/(n\epsilon)$ for
$n\in\bbZ$ meet $r=-\epsilon$ and $r=0$ at $\phi=\pi/2$.  If we
reflect the metric in $r=-\epsilon/2$,
\[\lambda(r)=\left\{\begin{array}{cc}
\cos^2(-r-\epsilon)&r\leq-\epsilon\\
1&-\epsilon<r<0\\\cos^2(r)&0\leq r \end{array}\right.\] then we have a
constructed an infinite family of geodesics with arbitrarily large
index that run from the boundary to the boundary.  These lift to
compact minimal surfaces in $S^3$.  The metric on $S^3$ is obtained
from the round metric by cutting $S^3$ along a totally geodesic $S^2$
and gluing in the product $[-\epsilon,0]\times S^2$ of the Euclidean
metric times the round metric.  The surfaces have arbitrarily large
index since by \cite{HLaMin} Jacobi fields lift to Jacobi fields
upstairs and we can use a generalisation of the Morse index theorem by
Simons \cite{SimMin} to get a lower bound for the index.

\section{Smooth perturbed metric.}
\begin{figure}[ht]
\begin{center}
\scalebox{0.5}{\includegraphics[height=8cm]{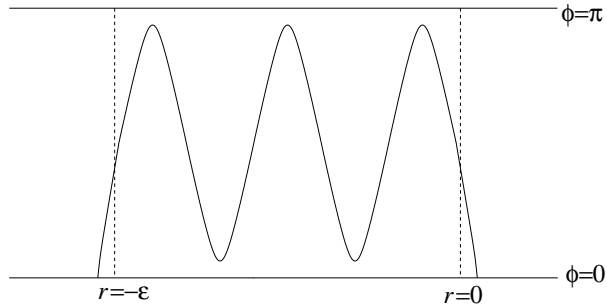}}
\end{center}
\caption{A geodesic meets the boundary twice.}
\label{fig:boundary}
\end{figure}

In this section we construct smooth metrics with similar properties to
the $C^1$ example in the previous section.  The main differences are
that we do not have an explicit formula for the geodesics and there is
no focal point.  We call a geodesic that meets the boundary a {\em
boundary geodesic.}  For particular smooth functions $\lambda(r)$
defined on a neighbourhood of $r=0$ the following properties hold:

\begin{enumerate}
    \item[(i)] Any boundary geodesic meets the boundary with known first
    and second derivatives.

    \item[(ii)] If we squeeze at a faster rate (a condition on
    $\lambda'(r)$) than for the $C^1$ example (\ref{eq:c1}), then the
    geodesics (\ref{eq:explicit}) form barriers for the boundary
    geodesics.
    
    \item[(iii)] There exist geodesics meeting all points of the boundary.

    \item[(iv)] The boundary geodesics meet the line $r=0$ at arbitrarily
    small angles.

    \item[(v)] There exists a countable set of boundary geodesics with
    arbitrarily high index that must meet the boundary again as in
    Figure~\ref{fig:boundary}.

    \item[(vi)] The Ricci curvature of the metric on the three-sphere can be
    chosen to be non-negative.
\end{enumerate}

In the remainder of this section we give the proofs of these
properties.  The precise meaning of barrier in (ii) is given in
Lemma~\ref{th:barrier}.  Proposition~\ref{th:exist} contains the proof
of (iii).  To prove (iv) we take geodesics meeting the boundary at
arbitrarily small $r>0$ and use (ii) to show these geodesics meet the
line $r=0$.  A second derivative estimate shows that the geodesics
remain close to vertical there.  As the geodesics become more
vertical, their period along $r<0$ becomes smaller and the index
becomes larger.  For (v) we prove that a countable set meet the lines
$r=-\epsilon$ and $r=0$ at the same value of $\phi$ and same angle
(with opposite sign.)  These pull back to embedded minimal surfaces in
$S^3$ of arbitrarily high index.

\begin{lemma}  \label{th:basic}
    There are geodesics meeting the boundary at $r>0$.  They meet the
    boundary orthogonally, with second derivative
    $d^2r/d\phi^2|_{\phi=0}=\lambda'(r)/2$ and each is the unique
    geodesic to meet its boundary point.
\end{lemma}
\begin{proof}
    The shortest path from any point to the boundary is a geodesic. 
    If $r>0$ then it must meet the boundary at $r>0$ since any
    geodesic that crosses $r=0$ cannot meet the boundary on $r<0$. 
    For $\phi$ small enough, there is a unique shortest path to the
    boundary by the uniqueness of a minimal disk upstairs with
    boundary given by a small circle.  A geodesic meets the boundary
    orthogonally since upstairs the circle action splits the tangent
    space of the fixed point into the invariant tangent plane of the
    minimal surface and its orthogonal complement pointing in the
    direction of the fixed point set.  In other words $dr/d\phi=0$ at
    the boundary.  The limit $\lim_{\phi\rightarrow 0}d^2r/d\phi^2$
    exists by smoothness of the minimal surface upstairs.  It is equal
    to $\lim_{\phi\rightarrow 0}(dr/d\phi)/\sin\phi$ and we can
    calculate it using (\ref{eq:implicit}): 
    \[ d^2r/d\phi^2|_{\phi=0}=-d^2r/d\phi^2|_{\phi=0}+\lambda'(r).\]
    Any point on the boundary is met by at most one geodesic: if not
    then two such geodesics would lift to tangent minimal surfaces
    violating the maximum principle.
\end{proof}

In the following lemma we prove that the geodesics from the $C^1$
example in the last section act as barriers for the geodesics of
metrics satisfying (\ref{eq:barrier}).  We refer to the former
geodesics as leaves of a foliation, to reduce
confusion when talking about two sets of geodesics.  The argument is
quite easy.  We prove that a geodesic of a metric that squeezes fast
enough can only be parallel from the left to a leaf of the foliation. 
A geodesic that moves away from the left of a leaf cannot return to
that leaf since it would have to meet another leaf at a tangent from
the right.  Here left refers to the $(r,\phi)$ cartesian plane with
$r$ on the horizontal axis.  The point $(r_1,\phi)$ lies to the left
of $(r_2,\phi)$ when $r_1<r_2$.
\begin{lemma}  \label{th:barrier}
    If
    \begin{equation}  \label{eq:barrier}
	\lambda'(r)< -\sin 2r
    \end{equation}
     for $r>0$ then a geodesic that meets the boundary at
     $(r,\phi)=(r_0,0)$ does not meet the curve $\cos\phi=\tan r/\tan
     r_0$ in $r\in[0,r_0)$.
\end{lemma}
\begin{proof}
    Along the leaf $\cos\phi=\tan r/\tan r_0$, at $(r,\phi)=(r_0,0)$
    we have $dr/d\phi=0$ and $d^2r/d\phi^2=-\cos r_0\sin r_0$.  By
    Lemma~\ref{th:basic}, a boundary geodesic meets $(r,\phi)=(r_0,0)$
    with $dr/d\phi=0$ and $d^2r/d\phi^2<-\cos r_0\sin r_0$.  Thus,
    near $\phi=0$ the boundary geodesic lies to the left of the leaf
    $\cos\phi=\tan r/\tan r_0$.  The leaves $\cos\phi=\tan
    r/\tan\kappa$ foliate the surface.  If the boundary geodesic were
    ever to meet the leaf $\cos\phi=\tan r/\tan r_0$ again then at
    some point it would be tangent to the foliation, from the right of
    a leaf.
    
    If $\lambda'(r)< -\sin 2r$ then $\lambda(r)<\cos^2r$ for $r>0$
    since $\lambda(0)=1$.  Along the boundary geodesic with direction
    $dr/d\phi=-\sin\phi\cos^2r\tan\kappa$ (the tangent direction to
    the foliation), from (\ref{eq:implicit}) we have
    \begin{eqnarray*}
	\frac{d^2r}{d\phi^2}&=&
	\frac{\cos^4r}{\lambda(r)}(\tan^2\kappa-\tan^2r)
	(\sin r\cos r +\frac{3}{2}\lambda'(r))+\sin r\cos r+\lambda'(r)\\
	&<&\cos^2r(\tan^2\kappa-\tan^2r)(-2\sin r\cos r)-\sin r\cos r
    \end{eqnarray*}
    where we have replaced $\lambda'(r)$ with $-\sin 2r$ since its
    coefficient is positive, and $1/\lambda(r)$ with $1/\cos^2r$ since
    its coefficient is negative.  Compare this with the second
    derivative of the foliation $r_{\rm fol}(\phi)$
    \[\frac{d^2r_{\rm fol}}{d\phi^2}=
    \sin r\cos r(1-2(1+\tan^2\kappa)\cos^2r)\]
    Then 
    \[\frac{d^2r_{\rm fol}}{d\phi^2}-\frac{d^2r}{d\phi^2}>0.\]
    
    Thus, if the geodesic is tangent to the foliation somewhere, the 
    second derivative of the foliation is greater than that of the 
    geodesic so the geodesic lies to the left of the leaf which 
    contradicts the geometry described above, that it should lie to 
    the right of the leaf.
\end{proof}

An implicit form of the geodesic equations is given in
(\ref{eq:implicit}).  In the following lemma we use this to show that
geodesics that meet the boundary at $\phi=0$ near $r=0$ remain almost
vertical.

\begin{lemma}  \label{th:monotone}
    If a geodesic meets the boundary at $\phi=0$, then along the
    geodesic, whilst $r\geq 0$, $dr/d\phi\leq 0$ and
    \begin{equation}  \label{eq:vertical}
	d^2r/d\phi^2>\lambda'(r).
    \end{equation} 
\end{lemma}
\begin{proof}
    By Lemma~\ref{th:barrier}, if a geodesic meets the boundary at
    $\phi=0$, then along the geodesic $\cos\phi/\sin\phi>0$.  This
    fact, together with $\lambda'(r)<0$ shows that the equation
    (\ref{eq:implicit}) expresses $d^2r/d\phi^2$ as a cubic polynomial
    in $dr/d\phi$ with negative coefficients.  Any such cubic
    polynomial $f=ax^3+bx^2+cx+d$, $a,b,c,d$ negative, satisfies
    $f(x)>d$ for $x<0$.  In our case, this yields
    $d^2r/d\phi^2>\lambda'(r)$ when $dr/d\phi<0$.  We can conclude
    that $dr/d\phi\leq 0$ for $r\geq 0$ as follows.  If $dr/d\phi>0$
    somewhere then as $r\rightarrow 0$, either $dr/d\phi=0$ somewhere
    or $dr/d\phi=\infty$ somewhere.  In the first case, $dr/d\phi=0$
    implies that $d^2r/d\phi^2=\lambda'(r)<0$ so immediately
    afterwards, $dr/d\phi<0$.  In the second case, the lower bound
    $d^2r/d\phi^2>\lambda'(r)$ prevents $|dr/d\phi|$ from becoming too
    big, by the mean value theorem (described in detail in the proof
    of the next Lemma.)
\end{proof}

\begin{lemma}  \label{th:convex}
    A boundary geodesic is convex in a neighbourhood of the boundary.
\end{lemma}
\begin{proof}
    The statement of the lemma is equivalent to the following: for
    small enough $\phi$ if a geodesic satisfies $d^2r/d\phi^2>0$
    somewhere then it doesn't meet the boundary at $\phi=0$.  From 
    this we deduce that a boundary geodesic satisfies 
    $d^2r/d\phi^2\leq 0$ near the boundary.
    
    As in the proof of Lemma~\ref{th:monotone} we use the cubic
    polynomial (\ref{eq:implicit}).  When $\phi$ is small enough the
    cubic polynomial is strictly decreasing.  To see this,
    differentiate (\ref{eq:implicit}) with respect to $dr/d\phi$. 
    This is quadratic in $dr/d\phi$ with discriminant
    \begin{equation}  \label{eq:discriminant}
	\Delta=\frac{3}{\lambda(r)}\left(
	\frac{3\lambda'(r)^2}{4\lambda(r)}
	-\frac{\cos^2\phi}{\sin^2\phi}\right).
    \end{equation}
    We require $\lambda'(r)=0$ when $\lambda(r)=0$ so
    $\lambda'(r)^2/\lambda(r)$ is bounded since its limit as
    $\lambda(r)$ vanishes is $2\lambda''(r)$.  Thus, for small enough
    $\phi$, (\ref{eq:discriminant}) is negative so the polynomial
    (\ref{eq:implicit}) is monotone.
    
    The unique zero of (\ref{eq:implicit}) is negative.  If
    $d^2r/d\phi^2>0$ somewhere then $dr/d\phi$ becomes more negative
    as $\phi$ decreases, and $d^2r/d\phi^2$ remains positive.  In
    particular, as $\phi$ decreases so does $dr/d\phi$.  But if the
    geodesic were to meet the boundary whilst $dr/d\phi\leq 0$ it
    would not meet with $dr/d\phi=0$ contradicting
    Lemma~\ref{th:basic}.  Thus it must repel from the boundary so
    $dr/d\phi>0$ and it cannot meet $\phi=0$ by
    Lemma~\ref{th:monotone}.
   
\end{proof}

\begin{prop}  \label{th:exist}
    There exist geodesics meeting all points of the boundary.
\end{prop}
\begin{proof}
    Take a sequence of points converging to $(r,\phi)=(r_0,0)$ and
    consider the shortest paths from each of these to the boundary. 
    By Lemma~\ref{th:basic}, when $\phi$ is small enough the shortest
    paths are unique.  The remainder of the proof shows that this set
    of geodesics converges to a geodesic that meets the boundary at
    $(r_0,0)$.
    
    The sequence of geodesics $r(\phi)$ and the sequence of
    derivatives $dr/d\phi$ are equicontinuous families so by
    Arzela-Ascoli converge uniformly to a geodesic.  Equicontinuity
    follows from boundedness and Lemmas~\ref{th:monotone} and
    \ref{th:convex} which supply a uniform bound on the second
    derivatives $d^2r/d\phi^2$ of the sequence of geodesics.  
    
    It remains to show that the limiting geodesic meets the boundary
    at $(r_0,0)$.  We will show that the target of each shortest
    geodesic lies close to $(r_0,0)$.
    
    In Lemma~\ref{th:convex} it was shown that when $\phi$ is small
    enough the polynomial (\ref{eq:implicit}) has a unique zero.  We
    can estimate it as follows.  Let $\delta$ be any small positive
    number.  Express (\ref{eq:implicit}) as
    \[\frac{d^2r}{d\phi^2}=\frac{1}{\lambda(r)}\left(-\frac{\cos\phi}
    {\sin\phi}\frac{dr}{d\phi}+\frac{3}{2}\lambda'(r)\right)
    \left(\frac{dr}{d\phi}\right)^2
    -\frac{\cos\phi}{\sin\phi}\frac{dr}{d\phi}+\lambda'(r)\] so if
    $dr/d\phi=-\delta$ then for small enough $\phi$ this expression is
    positive.  Since (\ref{eq:implicit}) is negative for $dr/d\phi=0$,
    the unique zero lies in the interval $(-\delta,0)$.
    
    From a point, say $(r_0,\phi)$ for $\phi$ chosen small enough as
    in the previous paragraph, the closest point on the boundary
    satisfies $r>r_0$ since $dr/d\phi<0$.  It also lies under any
    straight line $r/\phi=c<-\delta$ since the tangent direction to
    the shortest path must satisfy $dr/d\phi>-\delta$ so that
    $d^2r/d\phi^2<0$ and by the convexity proven in
    Lemma~\ref{th:convex}, the geodesic lies under its tangent line.
    
    Thus, as $\phi\rightarrow 0$, the target points on the boundary 
    converge to $(r_0,0)$. 
\end{proof}

\begin{lemma}   \label{th:small}
    The geodesics that meet the boundary meet the line $r=0$ at
    arbitrarily small angles.
\end{lemma}
\begin{proof}
    We have a geodesic running from $(0,\phi_0)$ to $(r_0,0)$
    satisfying $dr/d\phi\leq 0$ and $d^2r/d\phi^2>\lambda'(r)$ by
    Lemma~\ref{th:monotone}.  Put $dr/d\phi=\alpha$ at $r=0$ and
    recall that $dr/d\phi=0$ at $r=r_0$.  By the mean value theorem,
    there exists $0<r_1<r_0$ with $d^2r/d\phi^2=\alpha/\phi_0$ and
    since $\phi_0<\pi/2$ and $\alpha<0$ we get
    $d^2r/d\phi^2<2\alpha/\pi$ so by (\ref{eq:vertical}),
    \[ \lambda'(r)\pi/2<\alpha<0.\]  When $r_0$ is small, $\lambda'(r)$ is
    close to $0$ for $0<r<r_0$ so $\alpha$ is arbitrarily small.
    
\end{proof}

Reflect the metric in $r=-\epsilon/2$ and note that
$\lambda(r)=\lambda(-\epsilon-r)$ is still a smooth function.  Take
the geodesic that meets the boundary at $(r,\phi)=(r_0,0)$.  By
Lemmas~\ref{th:barrier} and \ref{th:small}, it meets the line $r=0$ at
a value $\phi_0<\pi/2$ with $dr/d\phi=\delta<0$.  If it also meets the
point $(r,\phi)=(-\epsilon,\phi_0)$ with $dr/d\phi>0$ then $dr/d\phi$
is necessarily $-\delta$ and the geodesic meets the boundary also at
$(-\epsilon-r_0,0)$.  The next lemma shows this occurs infinitely
often.
\begin{lemma}
    There exists a countable set of geodesics with arbitrarily
    high index that meet the boundary twice.
\end{lemma}
\begin{proof}
    Fix $\epsilon>0$.  When a geodesic meets $r=0$ at $\phi=\phi_0$
    with $dr/d\phi=\alpha<0$, then if $\alpha$ is small enough there
    are finitely many values $r_i$, $-\epsilon<r_i<0$ such that the
    geodesic passes through $(r_i,\phi_0)$ with $dr/d\phi=-\alpha$. 
    Order these points with $0>r_1>r_2>\ldots>r_k\geq-\epsilon$.  The
    difference $r_i-r_{i+1}$ gives the period and is constant.  The
    quantity $(r_k+\epsilon)/(r_i-r_{i+1})\in[0,1)$ varies
    continuously with the point $r_0$, where the geodesic meets the
    boundary.  For $r_0$ small enough, the angle between the geodesic
    and $r=0$ is small and so is the period.  Thus we can choose $r_0$
    small enough that there are more than $k$ values $r_i$ defined
    above.  By the intermediate value theorem, there is a value of
    $r_0$ in between such that the quantity
    $(r_k+\epsilon)/(r_i-r_{i+1})=0$.  This geodesic meets the
    boundary also at $(-\epsilon-r_0,0)$.  This construction produces
    such a geodesic for arbitrarily small $r_0$ and hence arbitrarily
    high index.
\end{proof}

\begin{lemma}  \label{th:ricci}
    If $\lambda$ satisfies (\ref{eq:barrier}) and $\lambda''(r)\leq
    2\sin^2r$ for $r>0$ then the Ricci curvature of the metric on the
    three-sphere is non-negative.
\end{lemma}
\begin{proof}
    The Ricci curvature of the metric is
    \begin{equation}
	\left(\begin{array}{ccc}
	\frac{1}{2}\lambda'(r)^2/\lambda(r)^2-\lambda''(r)/\lambda(r)&
	0&0\\0&(1-\frac{1}{2}\lambda''(r))/\lambda(r)&0\\0&0&
	(1-\frac{1}{2}\lambda''(r))/\lambda(r)\end{array}\right).
    \end{equation}
    In a neighbourhood of $\lambda(r)\equiv 1$ the first term of the
    Ricci curvature vanishes and the other two terms are 1.  For $r>0$
    (and $r<-\epsilon$), if $\lambda''(r)<0$ then clearly each term is
    positive.  Thus, assume $0\leq\lambda''(r)\leq 2\sin^2r$. 
    Together with $\lambda(r)<\cos^2r$ we get
    $-\lambda''(r)\lambda(r)>-2\sin^2r\cos^2r$.  Since
    $\lambda'(r)^2>\sin^22r$, then $(1/2)
    \lambda'(r)^2-\lambda''(r)\lambda(r)>0$ and the first entry of the
    Ricci curvature is positive.  The other two entries are bounded
    below by $\cos^2r/\lambda(r)>0$.
\end{proof}

Since the geodesics that meet the boundary arbitrarily close to $r=0$
are the important ones, we can take a part of the surface $r>r_0$ that
is untouched by the family of geodesics.  This corresponds to a ball
in the the three-sphere.  We can cut out this ball and glue in a
three-manifold so the minimal surfaces live in any three-manifold.

In the proof of the previous lemma it is shown that the Ricci
curvature is positive when $\lambda''(r)>0$.  Thus, in the interval
$-\epsilon\leq r\leq 0$ we can increase the sizes of the two-spheres
to get positive Ricci curvature.  That is, there is an arbitrarily
small deformation of $\lambda(r)$ in a neighbourhood of $-\epsilon\leq
r\leq 0$ with $\lambda''(r)>0$.  As mentioned in the introduction,
this bounds the area, and hence index, of each embedded minimal
sphere.

\end{document}